\documentclass[12pt]{article}
\usepackage[utf8]{inputenc}
\usepackage{amssymb,amsmath,amsfonts,eucal,mathrsfs,amsthm} 
\usepackage{hyperref}
\usepackage{empheq}
\newtheorem{theorem}{Theorem}
\newtheorem{proposition}[theorem]{Proposition}

\newtheorem{remark}[theorem]{Remark}

\theoremstyle{definition}
\newcommand{\R}{\mathbb{R}}

\newcommand{\Sf}{\mathbb{S}}
\newcommand{\C}{\mathbb{C}}

\newcommand{\grad}{\mbox{grad\,}}

\newcommand{\nab}{\tilde\nabla}

\newcommand{\End}{\mbox{End}}

\def\<{{\langle}}
\def\>{{\rangle}}

\def\T{{\cal T}}

\def\d{\partial}

\def\p{\partial}
\def\be{\begin{equation} }
\def\ee{\end{equation} }

\def\proof{\noindent{\it Proof:  }}
\def\qed{\ifhmode\unskip\nobreak\fi\ifmmode\ifinner
\else\hskip5 pt \fi\fi\hbox{\hskip5 pt \vrule width4 pt
height6 pt  depth1.5 pt \hskip 1pt }}
\makeatletter
\newcommand{\subjclass}[2][]{\let\@oldtitle\@title
\gdef\@title{\@oldtitle\footnotetext{#1 
\emph{Mathematics Subject Classification:} #2}}}
\newcommand{\keywords}[1]{\let\@@oldtitle\@title
\gdef\@title{\@@oldtitle\footnotetext
{\emph{Key words and phrases.} #1.}}}
\makeatother
\begin{document}

\date{}
\title{The infinitesimal deformations of hypersurfaces\\ 
that preserve the Gauss map}
\author{M. Dajczer and M. I. Jimenez}
\maketitle

\begin{abstract} Classifying the nonflat hypersurfaces in 
Euclidean space $f\colon M^n\to\R^{n+1}$ that locally admit  
smooth infinitesimal deformations that preserve the Gauss map 
infinitesimally was a problem only considered by Schouten 
\cite{Sc} in 1928. He found two conditions that are necessary 
and sufficient, with the first one  being the minimality of the 
submanifold. The second is a technical condition that does 
not clarify much about the geometric nature of the hypersurface. 
In that respect, the parametric solution of the problem given 
in this note yields that the submanifold has to be Kaehler.
\end{abstract}

The local study of the smooth deformations of the submanifolds
of the Euclidean space comprises several types, among them 
isometric, conformal and affine variations, as well as their 
somewhat less demanding infinitesimal versions. The case of 
surfaces was already a hot topic during the $19^{th}$ century, 
whereas work for higher dimensions can be traced back to the 
final decade of that century. A rather large amount of papers 
on these and even more general types of deformations have 
since been  published.  

Related to this paper, we recall that a necessary 
condition for a hypersurface $M^n$ in $\R^{n+1}$, $n\geq 3$, 
without flat points to admit a non-trivial infinitesimal variation 
is to have at any point precisely two nonzero principal curvatures. 
This rigidity result is already in the book by Cesàro \cite{Ce} 
published in 1896, but it seems that it was first obtained by 
Killing in 1885. 

A parametric local classification of the Euclidean hypersurfaces 
$f\colon M^n\to\R^{n+1}$, $n\geq 3$, that admit smooth 
isometric deformations was obtained by Sbrana \cite{Sb1} in 1909.
His achievement was the conclusion of a long quest for such 
a classification. Among several partial results there is 
one due to Bianchi \cite{Bi} for dimension $n=3$. 
As an aside, it seems that Bianchi had been Sbrana's adviser. 
Around the same time Sbrana \cite{Sb} obtained the local 
classification of the nonflat hypersurfaces that admit 
infinitesimal isometric deformations. These deformations are 
the ones that preserve lengths just ``up to the first order". 
At the very end of his paper, Sbrana showed that these 
hypersurfaces form a \emph{much larger} class than the one 
of the isometrically deformable hypersurfaces.

Still in the isometric case, in 1916 Cartan \cite{Ca1} obtained 
a classification equivalent to Sbrana's by means of his own methods. 
Then, the following year he gave in \cite{Ca2} a parametric 
classification for the more general conformal case. For modern 
versions of all these classical results for hypersurfaces as 
well as further developments we refer to \cite{DFT}, \cite{DJ1}, 
\cite{DJ2}, \cite{DJV}, \cite{DT} and \cite{DV}. 
\vspace{1ex}

The local basic question: \emph{To what extent a Euclidean 
submanifold in arbitrary codimension $f\colon M^n\to\R^{n+p}$, 
$n\geq 2$, is determined by its Gauss map into the Grassmannian 
$G_{np}$} was answered in \cite{DG1}, with the later complement 
given in \cite{DR}. Of course, for $\R^3$ there was 
the very classical result that a nonflat surface has to be 
minimal and that any other isometric immersion with the same Gauss 
map belongs to the associated one-parameter family given by 
the Weierstrass parametrization. For instance, see page 383 in 
Darboux \cite{Dar}. In the case of hypersurfaces 
$f\colon M^n\to\R^{n+1}$, $n\geq 3$, we have from part $(iii)$
of Corollary $4.12$ in \cite{DG1} that, besides the minimality 
of the immersion, the manifold has to be Kaehler up to
a possible cylindrical Euclidean factor. Moreover, the set 
of all its isometric deformations preserving the Gauss map is
the one-parameter associated family discussed with details in 
Section \ref{Kaehler}. See also the more general result given 
in \cite{DG2}. 

This paper is devoted to the study of the hypersurfaces 
$f\colon M^n\to\R^{n+1}$, $n\geq 2$, that admit infinitesimal 
isometric deformations that infinitesimally preserve the Gauss 
map. Although we were not able to find in the literature 
any previous attempt  to give a classification in the isometric 
case other than the one aforementioned, it was for us quite 
unexpected to discover that there is one for the infinitesimal 
case. In fact, in terms of a rather difficult notation Schouten 
\cite{Sc} gave in 1928 a necessary and sufficient condition.
His main achievement is that we must be in the presence of a 
minimal hypersurface. Then, there is a second condition 
somehow involving the second fundamental form of the submanifold. 
Of course, Schouten's result does not fully describe 
the hypersurfaces involved.

The following result gives a complete local parametric 
classification of the Euclidean hypersurfaces that admit 
infinitesimal isometric deformations which infinitesimally 
preserve the Gauss map. Although, as aforementioned, while 
the class of hypersurfaces in the infinitesimally deformable 
class is much richer than in the isometrically deformable 
class, it turns out that that under the Gauss map condition
we are reduced to the same family of hypersurfaces, namely, 
the minimal ones that, up to a possible Euclidean cylindrical 
factor, are Kaehler. Since the presence of the former does 
not really significantly change the situation, that possibility 
is not considered in the following statement.

\begin{theorem}\label{main} Let $f\colon M^n\to\R^{n+1}$, $n\geq 2$, 
be an isometric immersion that admits a nowhere trivial infinitesimal 
variation that infinitesimally preserves the Gauss map. Assume that 
$M^n$ is free of flat points and that $f$ is not part of a 
Euclidean cylinder when restricted to any open subset of $M^n$.
Then $f$ is a minimal immersion of a Kaehler manifold. In addition, 
if $M^n$ is simply connected, then the variational vector 
field associated to the infinitesimal variation is, up to constants, 
the conjugate to $f$ in its associate family.
\end{theorem}

The assumptions that a hypersurface has no local Euclidean factor, 
as well as that the infinitesimal variation is nowhere trivial, are 
explained in the next section. Two local parametric descriptions 
of the minimal Kaehler hypersurfaces are provided in Section 2, among 
them a Weierstrass type representation that extends the classical one 
for minimal surfaces.

\section{Infinitesimal isometric variations}

On this section, we recall from \cite{DJ1} and \cite{DT} several 
concepts as well as some  of the basic facts on infinitesimal 
deformations of Euclidean hypersurfaces.
\vspace{1ex}

Let $F\colon I\times M^n\to\R^{n+1}$ denote a smooth 
variation of a given isometric immersion $f\colon M^n\to\R^{n+1}$ 
of a Riemannian manifold $M^n$, $n\geq 2$, such that 
$0\in I\subset\R$ is an interval and $F(0,x)=f(x)$. We 
assume that $f_t=F(t,\cdot)\colon M^n\to\R^{n+1}$ is an immersion 
for any  $t\in I$.  The variational vector field of the variation 
$F$ of $f$ is the section $\T\in\Gamma(f^*T\R^{n+1})$ defined as 
$$
\T=F_*\d/\d t|_{t=0}=\nab_{\d/\d t}f_t|_{t=0},
$$
where $\nab$ denotes the Euclidean connection.

We assume that $F\colon I\times M^n\to\R^{n+1}$ is an 
\emph{infinitesimal variation}. This means that the metrics 
$g_t$, $t\in I$, induced by $f_t$ satisfy 
$\partial/\partial t|_{t=0}g_t=0$, that is, 
\be\label{varcond2}
\frac{\d}{\d t}|_{t=0}\<f_{t*}X,f_{t*}Y\>=0\;\;
\mbox{for any}\;\; X,Y\in\mathfrak{X}(M). 
\ee
That is what it is meant by saying that the variation preserves 
lengths ``up to the first order". Equivalent to \eqref{varcond2} 
is that the variational vector field $\T\in\Gamma(f^*T\R^{n+1})$ 
of $F$ satisfies the condition
\be\label{infben}
\<\T_*X,f_*Y\>+\<f_*X,\T_*Y\>=0\;\;
\mbox{for any}\;\; X,Y\in\mathfrak{X}(M).
\ee

A vector field $\T\in \Gamma(f^*T\R^{n+1})$ satisfying \eqref{infben} 
is called an \emph{infinitesimal bending} of $f$. It is known 
already from classical differential geometry that the convenient 
approach to study infinitesimal variations is to focus on the 
infinitesimal bending determined by the variation.

Associated to a given infinitesimal bending $\T$ of 
$f\colon M^n\to\R^{n+1}$ we have the infinitesimal variation 
$F\colon\R\times M^n\to\R^{n+1}$ given by
\be\label{unique0}
F(t,x)=f(x)+t\T(x)
\ee
which has $\T$ as variational vector field. 
By no means this infinitesimal variation is unique with 
this property, although it may be seen as the simplest one. 
In fact, new infinitesimal variations with variational vector 
field $\T$ are obtained by adding to \eqref{unique0} terms of the 
type $t^k\delta$, $k\geq 2$, where $\delta\in\Gamma(f^*T\R^{n+1})$.
\vspace{1ex}

An infinitesimal bending is called a \emph{trivial} infinitesimal 
bending if it is the restriction to the submanifold of a Killing 
vector field of the ambient space. That is, it is the variational 
vector field of a variation by Euclidean motions.
Hence, there is a skew-symmetric linear endomorphism 
${\cal D}$ of $\R^{n+1}$ and a vector $w\in\R^{n+1}$ such that 
${\cal T}={\cal D}f+w$. Conversely, given a trivial infinitesimal 
bending we have the trivial isometric variation of $f$ given by 
$F(t,x)=e^{t{\cal D}}f(x)+tw$.

Let $A(t)$, $t\in I$, with $A(0)=A$ denote the shape operators 
associated to the smooth family of Gauss maps $N_t$ of $f_t$ 
where $N_0=N$. 
Then let $B\in\Gamma(\End(TM))$ be the symmetric Codazzi tensor 
defined by $B=\p/\p t\vert_{t=0}A(t)$. 

The Gauss formula for each $f_t$  states 
\be\label{Gaussform}
\nab_Xf_{t*}Y=f_{t*}\nabla^t_XY+\<f_{t*}A_tX,f_{t*}Y\>N_t
\;\;\mbox{for any}\;\;X,Y\in\mathfrak{X}(M),
\ee
where $\nabla^t$ is the Levi-Civita connection of the metric 
induced by $f_t$. Notice that
$$
\partial/\partial t|_{t=0}\nab_Xf_{t*}Y
=\partial/\partial t|_{t=0}\nab_X\nab_Yf_t=\nab_X\T_*Y.
$$
Since $F$ is an infinitesimal variation then 
$\partial/\partial t|_{t=0}\nabla^t_XY=0$. 
Thus the derivative with respect to $t$ at $t=0$ of 
the right hand side of \eqref{Gaussform} gives
$$
\partial/\partial t|_{t=0}(f_{t*}\nabla^t_XY
+\<f_{t*}A_tX,f_{t*}Y\>N_t)
=\T_*\nabla_XY+\<BX,Y\>N+\<AX,Y\>\partial/\partial t|_{t=0}N_t,
$$
where $\nabla=\nabla^0$ is the Levi-Civita connection of the metric induced by $f$.
Hence, 
\be\label{B}
(\nab_X\T_*)Y=\<AX,Y\>\partial/\partial t|_{t=0}N_t+\<BX,Y\>N
\;\;\mbox{for any}\;\;X,Y\in\mathfrak{X}(M),
\ee
where $(\nab_X\T_*)Y=\nab_X\T_*Y-\T_*\nabla_XY$.
\vspace{1ex}

In addition, we have from equation $(2.21)$ in \cite{DJ1} that 
the tensors $A$ and $B$ satisfy 
\be\label{infGauss2}
BX\wedge AY-BY\wedge AX=0\;\;\mbox{for any}\;\;X,Y\in\mathfrak{X}(M).
\ee
This equation together with the Codazzi equation for $B$ form  
what is called the fundamental equations of an infinitesimal 
bending for a hypersurface.\vspace{1ex}

The following result is Proposition $2.10$ in \cite{DJ1}.

\begin{proposition}\label{trivial} The infinitesimal bending 
$\T$ is trivial if and only if $B=0$.
\end{proposition}

Let $N_t$ be the Gauss maps of $f_t$, $t\in I$.   
That the non-trivial infinitesimal variation 
$F$ of $f$ is \emph{infinitesimal with respect to the Gauss map}
means that $\partial/\partial t|_{t=0}N_t=0$. In other words,
$N_t$ coincides with the Gauss map $N$ of $f$ up to the first order. 
In this situation, the derivative with respect to 
$t$ of $\<N_t,f_{t*}X\>=0$ computed at $t=0$ yields
\be\label{normal}
\<N,\T_*X\>=0\;\;\text{for any}\;\;X\in\mathfrak{X}(M). 
\ee
Hence $\T_*X\in\Gamma(f_*TM)$ for any $X\in\mathfrak{X}(M)$, 
and thus determines a skew symmetric endomorphism of 
the tangent space also denoted by $\T_*\in\Gamma(\End(TM))$ 
for simplicity. 
\vspace{1ex}

Let $M^n$ be a Riemannian product $N^k\times V$ where $N^k$ 
is a simply connected manifold free of flat points
and $V\subset\R^{n-k}$ is an open subset. Let
$g\colon N^k\to\R^{k+1}$ be an isometric immersion and 
$f\colon M^n\to\R^{k+1}\oplus\R^{n-k}$ the cylinder
given by $f(y,z)=(g(y),z)$. Let $G$ be an infinitesimal
variation of $g$ and let $F$ be the variation of $f$
defined by $F(t,y,z)=(G(t,y),\Psi_t(z))$, where $\Psi_t$ is a
smooth one-parameter family of isometries of $\R^{n-k}$
with $\Psi_0=I$. Then $F$ is an infinitesimal variation
of $f$ by cylinders. The corresponding infinitesimal
bending is of the form
\be\label{bending}
\T=(\T_1,\T_2)
\ee
where $\T_1$ corresponds to $G$ and $\T_2$ is a Killing 
vector field of $V$.
The next result shows that this trivial 
situation always holds in the presence of a cylinder and 
explains why the existence of a local Euclidean cylindrical factor
has been excluded in the statement of Theorem \ref{main}.

\begin{proposition} Any non-trivial infinitesimal bending of 
$f$ as above is as in \eqref{bending} up to a trivial 
infinitesimal bending.
\end{proposition}

\proof Let $\T$ be a nontrivial infinitesimal bending of $f$ 
with associated tensor $B$. From Proposition \ref{trivial} 
we know that $B\neq 0$. Since $M^n$ is free of flat points, 
it follows from \eqref{infGauss2} that $\Delta=\ker A$ 
satisfies $\Delta\subset\ker B$ at any point of $M^n$. 
In particular, we have that $BS=0$ for 
any $S\in\Gamma(TV)$. Then the  Codazzi equation 
$(\nabla_SB)X=(\nabla_XB)S$
for $X\in\Gamma(\Delta^\perp)$ gives that $B$ is a parallel 
tensor along $V$. Recall that the second fundamental 
form of $f$ is given by the second fundamental 
form of $g$ on $\Delta^\perp$. Thus $B$ is a Codazzi tensor 
defined on  $N^k$ which, together with the second fundamental 
form of $g$, satisfies \eqref{infGauss2}. Hence, it follows 
from the Fundamental Theorem of infinitesimal bendings, namely, 
Theorem $2.11$ in \cite{DJ1}, that $B$ is the 
associated tensor of an infinitesimal bending $\T_1$ of $g$.
Finally, it follows from Proposition \ref{trivial} that 
$\T-(\T_1,0)$ is a trivial infinitesimal bending of $f$.\qed

\section{Euclidean Kaehler hypersurfaces}\label{Kaehler}

In this section, we first discuss a local parametrization of
the nonflat Euclidean Kaehler hypersurfaces in terms of the 
so called Gauss parametrization, in particular, when 
the submanifold is minimal. Then, we give a Weierstrass type 
representation for the hypersurfaces in the latter situation.
\vspace{1ex}

Let $f\colon M^n\to\R^{n+1}$, $n\geq 3$, be a Euclidean 
hypersurface of rank two, that is, with precisely two nonzero 
principal curvatures at each point. Such a hypersurface can be 
locally parametrized by means of the Gauss 
parametrization, first used by Sbrana in his two aforementioned 
papers, and given by \eqref{Gausspar} bellow. For additional 
details on this parametrization we refer to Theorem $7.18$ 
in \cite{DT}. 

Let $g\colon L^2\to\Sf^n$ be an isometric immersion of a 
surface into the unit sphere  with normal bundle $\Lambda$ 
and let $\gamma\in C^\infty(L^2)$ be an arbitrary function. 
Set $h=i\circ g$ where $i\colon\Sf^n\to\R^{n+1}$ denotes 
the inclusion map. Then the smooth map 
$\Psi\colon\Lambda\to\R^{n+1}$ given by
\be\label{Gausspar}
\Psi(x,w)=\gamma(x)h(x)+h_*\grad\gamma(x)+i_*w
\ee
is along the open subset of $\Lambda$ of regular points a 
parametrization of a hypersurface of rank two.  Conversely,
any hypersurface of rank two admits such a local parametrization.

According to Theorem $2.5$ in \cite{DG1} or Theorem $15.14$ in 
\cite{DT} the hypersurface is Kaehler if and only if the surface 
$g$ is pseudoholomorphic in the sense of Calabi \cite{Ca}. 
The latter means that there exists an orthogonal tensor $T$ in 
$\Lambda$ that is parallel with respect to the normal connection 
and satisfies $A_{T\xi}=J_0\circ A_\xi$ for all $\xi\in\Lambda$,
where $J_0$ is the almost complex structure on $L^2$. In 
particular, the surface $g$ is minimal.

The submanifold given by \eqref{Gausspar} is minimal if and 
only if $\Delta\gamma+2\gamma=0$ is satisfied. Notice that 
in this case, the hypersurface is real analytic. Moreover,  
the Gauss map of the hypersurface is $N(x,w)=g(x)$ and 
$\gamma=\<f,N\>$ is its support function. 

From Section $15.3.1$ in \cite{DT} we have that any minimal 
simply connected Kaehler hypersurface without flat points 
$f\colon M^n\to\R^{n+1}$  admits a one-parameter 
associated family of noncongruent isometric minimal immersions 
$f_\theta\colon M^n\to\R^{n+1}$, $\theta\in[0,\pi)$ with $f_0=f$. 
In fact, if $J_\theta=\cos\theta I + \sin\theta J$ then 
the \emph{associated family} is given by the line integral 
\be\label{assfam}
f_\theta(x)=\int_{x_0}^xf_*\circ J_\theta
\ee
where $x_0$ is any fixed point of $M^n$. Since \eqref{assfam} 
yields that ${f_\theta}_*=f_*\circ J_\theta$ then all the 
hypersurfaces in the family have the same Gauss map. Moreover, 
the hypersurface can be realized as the real part of its
\emph{holomorphic representative} $F\colon M^n\to\C^{n+1}$
given by $\sqrt{2}F(x)=(f(x),\bar{f}(x))
\in\R^{n+1}\oplus\R^{n+1}\cong\C^{n+1}$, where $\bar{f}
=f_{\pi/2}$ is called the \emph{conjugate} hypersurface
to $f$. Finally, we have that $A(\theta)=A\circ J_\theta$
is the shape operator of $f_\theta$. 
\vspace{1ex}

In the unpublished work of Hennes \cite{He} it was observed 
that the arguments used in \cite{DG3} can be adopted in order 
to produce the local Weierstrass type representation of all 
nonflat minimal Kaehler hypersurfaces discussed next.
\vspace{1ex}

Let $U\subset\C$ be a simply connected domain. Start with
a nonzero holomorphic function $\alpha_0\colon U\to\C$ and 
let $\phi_0=\int\alpha_0(z)dz$. Assuming that the functions 
$\alpha_r, \phi_r\colon U\to\C^{2r+1}$ have been defined for
some $0\leq r\leq n-1$ choose any nonzero holomorphic 
function $\mu_{r+1}\colon U\to\C$ and set
$$
\alpha_{r+1}=\mu_{r+1}
\begin{pmatrix}
\frac{1-\phi_r^2}{2}\\ i\frac{1+\phi_r^2}{2}\\ \phi_r\end{pmatrix}
\;\;\mbox{and}\;\;\phi_{r+1}=\int\alpha_{r+1}(z)dz,
$$
where $\phi_r^2=\phi_r\cdot\phi_r$ with respect to the standard 
symmetric inner product in $\C^{2r+1}$.

Let $W$ be an open subset in $\C^{n-1}$  containing the origin.  
Then let $F\colon U\times W\to\C^{n+1}$ be given by
$$
F(z,w_1,\ldots,w_{n-1})
=\sum_{j=0}^{n-1}\int b_j(z)\delta^{(j)}(z)dz 
+\sum_{j=1}^{n-1}w_j\delta^{(j-1)}(z)
$$
where $\delta=\alpha_n$, $\delta^{(j)}=d^j\delta/dz^j$ and 
$b_0,b_1,\ldots,b_{n-1}\colon U\to\C$ are holomorphic 
functions such that $b_{n-1}$ is never zero. Set $w_j=u_j+iv_j$,
$u=(u_1,\ldots,u_{n-1})$ and $v=(v_1,\ldots,v_{n-1})$.
Let $f(z,w_1,\ldots,w_{n-1})=\sqrt{2}Re(F(z,w_1,\ldots,w_{n-1}))$, 
that is,
\be\label{param}
\frac{1}{\sqrt{2}}f(z,u,v)
=Re\sum_{j=0}^{n-1}\int b_j(z)\delta^{(j)}(z)dz +\sum_{j=1}^{n-1}
\left(u_j Re\,\delta^{(j-1)}(z)-v_jIm\,\delta^{(j-1)}(z)\right).
\ee
Then $M^{2n}=(U\times W, f^*\<\,,\,\>)$ is a Kaehler manifold and 
$f\colon M^{2n}\to\R^{2n+1}$ is minimal. 
The conjugate minimal hypersurface of $f$ is 
$\bar{f}=\sqrt{2}Im(F)$. Hence $F$ is the 
holomorphic representative of $f$. Moreover, 
$f_\theta=\cos\theta f+\sin\theta \bar{f}$ is the full associated
family.   
Since $d f_\theta/d\theta|_{\theta=0}=\bar{f}$ thus the conjugate
hypersurface $\bar{f}=f_{\pi/2}$  in the associated family to $f$ 
is the variational vector field of the variation $f_\theta$.  

\begin{remark} {\em Two examples on case $M^4$ in $\R^5$
of \eqref{param} are computed explicitly in \cite{He}.
}\end{remark}

\section{The proof of Theorem \ref{main}}

We first prove that the hypersurface is minimal.
We have that the non-trivial infinitesimal variation
$F$ of $f$ is infinitesimal with respect to the Gauss map.
Then \eqref{B} yields
$$
(\nab_X\T_*)Y=\<BX,Y\>N\;\;\text{for any}\;\;X,Y\in\mathfrak{X}(M).
$$
We have seen that $\T_*\in\Gamma(\End(TM))$ is
a parallel tensor on $M^n$, this is, 
\be\label{Tpar}
(\nabla_X\T_*)Y=0\;\;\mbox{for any}\;\;X,Y\in\mathfrak{X}(M).
\ee
Then the above implies that $\<AX,\T_*Y\>=\<BX,Y\>$ 
for any $X,Y\in\mathfrak{X}(M)$. In other terms, we have
\be\label{BAT}
B=A\circ\T_*=-\T_*\circ A.
\ee

If the second fundamental form has at least rank three 
at any point, we have from the classical rigidity result 
discussed in the introduction that the infinitesimal bending 
is trivial; see Theorem $2$ in \cite{DR2} or Theorem $14.4$ 
in \cite{DT} . Hence, our assumptions imply that the 
hypersurface possesses constant rank $2$. 
It follows from \eqref{infGauss2} that the relative nullity 
vector subspace $\Delta=\ker A$ satisfies $\Delta\subset\ker B$ 
at any point. 
We have from \eqref{BAT} that $\T_*$ leaves $\Delta$ invariant 
and being $\T_*$ skew-symmetric also $\Delta^\perp$ 
is left invariant. Since $\Delta^\perp$ is two-dimensional, 
then we necessarily have that $\T_*|_{\Delta^\perp}$ is a 
multiple of a rotation of angle $\pi/2$, say 
\be\label{rot}
\T_*|_{\Delta^\perp}=cR_{\pi/2},
\ee
where $c\in\R$ from \eqref{Tpar}. Proposition \ref{trivial} 
gives that $B\neq0$ and hence $c\neq 0$.
Expressing $A$ and $\T_*|_{\Delta^\perp}$ in an 
orthonormal basis of principal directions, we have from
\eqref{BAT} that
$$
B|_{\Delta^\perp}=A|_{\Delta^\perp}\T_*|_{\Delta^\perp}
=\begin{bmatrix}
\lambda_1 & 0\\
0 & \lambda_2
\end{bmatrix}
\begin{bmatrix}
0 & -c\\
c & 0
\end{bmatrix}=
c\begin{bmatrix}
0 & -\lambda_1\\
\lambda_2 & 0
\end{bmatrix}.
$$
Then the symmetry of $B$ yields that $\lambda_1+\lambda_2=0$,
that is, $f$ is a minimal hypersurface and, in particular, 
real analytic. Notice that this concludes the proof for the 
surface case. 

We now prove that the $M^n$, $n\geq 3$, is Kaehler.
We have from \eqref{infben} that the vector subspaces 
$\ker \T_*$ and $(\ker\T_*)^\perp$ are $\T_*$-invariant. 
Then, on any open subset of $M^n$ where their dimensions 
are constant it follows from \eqref{Tpar} that they form 
totally geodesic distributions. Since 
$\T_*|_{\Delta^\perp}\neq 0$ and leaves  $\Delta^\perp$  
invariant then $\ker\T_*\subset\Delta$. Hence, if 
$\dim\ker\T_*=\ell>0$ is constant, we necessarily have that 
$f$ is locally a cylinder over an isometric immersion 
$g\colon N^{2k}\to\R^{2k+1}$ of a  Riemannian manifold 
$N^{2k}$ with $2k=n-\ell$. But this possibility has been  
ruled out by assumption.

Let $f_t\colon M^n\to\R^{n+1}$ be the variation of $f$ defined
by $f_t=f+t\T$ for $t\in\R$. Using \eqref{infben} we have 
$$
\|f_{t*}X\|^2=\|f_*X\|^2+t^2\|\T_*X\|^2.
$$
Thus $f_t$ and $f_{-t}$ for each $t$ induce the same 
Riemannian metric $g_t$ on $M^n$. Moreover, being $\T$ 
non-trivial it follows from either Proposition $14.3$ in 
\cite{DT} or Proposition $2.12$ in \cite{DJ1} that these 
two immersions are not congruent. In our situation, we 
have from \eqref{normal} that the unit normal vector field 
$N$ is in fact normal to $f_t$ for all $t$, i.e., 
the $f_t$, $t\in\R$, share the same Gauss map. 

As already observed in the introduction, it now follows 
from Corollary $4.12$ in \cite{DG1} that the Riemannian 
manifold $(M^n,g_t)$ is Kaehler, hence it admits an almost 
complex structure $J_t$ which is a parallel orthogonal tensor, 
and that the immersions $f_t$ and $f_{-t}$ are both minimal. 
Then, from Section \ref{Kaehler} it follows that the Gauss map 
$N_t=N$ of $f_t$ determines a pseudoholomorphic surface 
of $\Sf^n$. Since $f$ shares the same Gauss map with all 
the $f_t$, thus $(M^n,g_0)$ is a Kaehler manifold of 
dimension $n=2k$. 

Finally, we characterize the variational vector field.
Let $J$ be the almost complex structure of $M^{2k}$. Since 
the immersion $f$ is minimal we have from either Theorem $1.2$ 
in \cite{DR} or Theorem $15.7$ in \cite{DT} that $A\circ J=-J\circ A$. 
Therefore $\Delta$ and $\Delta^\perp$ are  $J$-invariant 
and thus $J|_{\Delta^\perp}=R_{\pi/2}$. 
Hence, it follows from \eqref{assfam} that the conjugate minimal 
immersion $\bar{f}=f_{\pi/2}$ satisfies 
$\bar{f}_*X=f_*JX$ for any $X\in\mathfrak{X}(M)$.
In particular, 
$$
\<f_*X,\bar{f}_*X\>=0\;\;\mbox{for any}\;\;X\in\mathfrak{X}(M)
$$
and thus $\bar{f}$ is an infinitesimal bending of $f$. 
Moreover, recall from Section \ref{Kaehler} that the 
associated family to $f$ has $\bar{f}$ as its variational 
vector field  and that all the immersions share the same 
Gauss map. Thus, in particular, \eqref{normal} holds for 
$\bar{f}$. 
It follows from \eqref{B} that its associated tensor $\bar{B}$ 
is the second fundamental form of $\bar{f}$, that is, 
$\bar{B}=A\circ J$.  On the other hand, from \eqref{BAT} 
and \eqref{rot} the tensor associated to $\T$ satisfies 
$B|_{\Delta^\perp}=cA|_{\Delta^\perp}\circ R_{\pi/2}$. 
Then Proposition \ref{trivial} yields that the infinitesimal 
bending $\T-c\bar{f}$ is trivial. Therefore, we have
$\T-c\bar{f}=\mathcal{D}f+w$
where $\mathcal{D}\in\End(\R^{2k+1})$ is skew-symmetric and 
$w\in\R^{2k+1}$. Since $\T_*$ and $c\bar{f}_*$ coincide on 
$\Delta^\perp$ then $\Delta^\perp\subset\ker\mathcal{D}$.
Moreover, we have from \eqref{normal} for $\T$ and $\bar{f}$ 
that $\<\mathcal{D}X,N\>=0$
for all $X\in\mathfrak{X}(M)$ and hence $N\in\ker\mathcal{D}$.
Therefore $\ker\mathcal{D}^\perp$, seen in $T_{f(x)}\R^{2k+1}$, 
lies in $f_*\Delta(x)$ at any point $x\in M^{2k}$. Being 
$\mathcal{D}$ constant we have that $f$ would necessarily 
be a cylinder with Euclidean factor given by $\ker\mathcal{D}^\perp$. 
Thus $\mathcal{D}$ vanishes and therefore $\T$ coincides with 
$c\bar{f}+w$ for some $c\in\R$ and $w\in\R^{2k+1}$.\vspace{2ex}\qed

\noindent \textbf{Funding:} Marcos Dajczer is partially 
supported by the grant PID2021-124157NB-I00 funded 
by MCIN/AEI/10.13039/501100011033/ ‘ERDF A way of making Europe’,
Spain, and is also supported by Comunidad Autónoma de la Región 
de Murcia, Spain, within the framework of the Regional Programme 
in Promotion of the Scientific and Technical Research (Action Plan 2022), 
by Fundación Séneca, Regional Agency of Science and Technology, 
REF, 21899/PI/22. Miguel I. Jimenez is supported by FAPESP with 
the grant 2022/05321-9.

\noindent Marcos Dajczer\\
Departamento de Matemáticas\\ 
Universidad de Murcia, Campus de Espinardo\\ 
E-30100 Espinardo, Murcia, Spain\\
e-mail: marcos@impa.br
\bigskip

\noindent Miguel Ibieta Jimenez\\
Instituto de Ciências Matemáticas e de Computação\\
Universidade de São Paulo\\
São Carlos\\
SP 13566-590-- Brazil\\
e-mail: mibieta@icmc.usp.br

\end{document}